\documentclass[11pt]{amsart}
\usepackage{amsmath,amssymb,mathrsfs,fontenc,moreverb} 
\usepackage{amsfonts}
\usepackage{array,delarray}
\usepackage[dvips]{graphicx} 
\usepackage[centredisplay,PostScript=dvips]{diagrams} 
     
\newcommand{\ass}{\mathrm{Ass}}
\newcommand{\codim}{\mathrm{codim}}
 
\newtheorem{theorem}{Theorem}[section]   
\newtheorem{corollary}[theorem]{Corollary}     

\theoremstyle{definition}

\theoremstyle{remark}
\newtheorem{example}[theorem]{Example}        

\numberwithin{equation}{section}  

\begin{document}

\title[Integral Closures of Cohen-Macaulay monomial ideals]  
{Integral Closures of Cohen-Macaulay monomial ideals} 
\author[Abdul Salam Jarrah]{Abdul Salam Jarrah}  
\address[Abdul Salam Jarrah]{Department of Mathematical Sciences\\
New Mexico State University\\
Las Cruces, NM 88003, USA}    
\email{ajarrah@nmsu.edu}     


\date{December 4, 2001}

\begin{abstract}
The purpose of this paper is to 
 present  a family of  Cohen-Macaulay monomial ideals such that their 
integral closures have embedded components and hence are not Cohen-Macaulay.  
\end{abstract}
\maketitle
\section{Introduction}
\noindent
Let $R$ be an arbitrary ring. 
Krull \cite{K} asked whether there exists a primary ideal $I$ 
such that the integral closure $\overline{I}$ of $I$ has embedded primes. 
Huneke \cite{H} gave a counterexample, his example  
in a 3-dimensional regular local ring of characteristic 2. 
No other examples have appeared in the literature. The long history
of Krull's question and the machinery that is required to 
understand Huneke's example suggest that there might not
be an easy counterexample. 
In fact,
it is easy to see that if $I$ is a $P$-primary monomial ideal in a polynomial 
ring, then the integral closure $\overline{I}$ of $I$ 
is a $P$-primary monomial ideal.
Thus, there is no monomial ideal that can serve as a
counterexample to Krull's question. 
Vasconcelos asked 
a relaxed version of Krull's question, namely whether there exists an
unmixed monomial ideal $I$ in a polynomial ring 
such that the integral closure $\overline{I}$ has embedded primes. 

In this paper, 
I present a family of Cohen-Macaulay monomial ideals such that
their integral closures are not Cohen-Macaulay, in particular,
their integral closures have embedded primes. Thus, each example in
this family is a counterexample to Vasconcelos's question.
Naturally, this family is not the only family of monomial
ideals with this property. Another Cohen-Macaulay monomial ideal 
is given at the end which does not belong to the family,
but whose integral closure has an embedded prime. 

Throughout this paper, $n$ and $t$ are  
positive integers and $n \geq 3$.

\section{Examples}
\noindent 
Let $R = k[x_1,\dots,x_n]=k[{\bf{\tt x}}]$ be a polynomial  
ring over an arbitrary field $k$. For any positive integer
$m$, let $[m] = \{1,\dots, m\}$.
For $i \in [n]$, let $e_i = (t,\dots,t,0,t,\dots,t) \in \mathbb{R}^n $,
with $0$ is in the $i$-th coordinate. Thus the corresponding monomial
with exponent $e_i$ is 
${\bf{\tt x}}^{e_i} = x_1^t \cdots x_{i-1}^t x_{i+1}^t \cdots x_n^t$.
Consider the monomial ideal
\[
I_{n,t} = \langle {\bf \tt x}^{e_1}, \dots, {\bf \tt x}^{e_n}\rangle.
\]
For $i, j \in [n]$ and $i\neq j$, 
it is easy to see  that  $I_{n,t} \subset \langle x_i^t, x_j^t \rangle$.
Thus 
$I_{n,t} \subseteq  \bigcap_{i=1}^{n-1} \bigcap_{j=i+1}^{n}\langle x_i^t, x_j^t \rangle$ and it is easy to see that equality holds. 
Moreover, without lose of generality, 
$\bigcap_{i=2}^{n-1} \bigcap_{j=i+1}^{n}\langle x_i^t, x_j^t \rangle \nsubseteq  \langle x_1^t, x_2^t \rangle$, since $x_3^t \cdots x_n^t \in \bigcap_{i=2}^{n-1} \bigcap_{j=i+1}^{n}\langle x_i^t, x_j^t \rangle$ but $x_3^t \cdots x_n^t \notin \langle x_1^t, x_2^t \rangle$. Therefore,
\begin{equation*}
I_{n,t} = \bigcap_{i=1}^{n-1} \bigcap_{j=i+1}^{n}\langle x_i^t, x_j^t \rangle
\end{equation*}
is the minimal primary decomposition of $I_{n,t}$.
In particular, $I_{n,t}$ has no embedded primes and $\codim(I_{n,t}) = 2$. On the other hand, 
the minimal free resolution of $R/{I_{n,t}}$ is
\[
0  \leftarrow  R/{I_{n,t}}
\leftarrow 
R
\xleftarrow{\begin{array}({ccc}) {\bf \tt x}^{e_1} &\cdots&  {\bf \tt x}^{e_n} \end{array}} 
R^n \xleftarrow{
\begin{array}({cccc}) 
-x_1^t & -x_1^t & \cdots & -x_1^t\\  
x_2^t & 0 & \cdots & 0\\ 
0 & x_3^t& \cdots & 0\\ 
\vdots & \vdots & \ddots & \vdots\\ 
0 & 0 &  \cdots & x_n^t
\end{array}} 
R^{n-1} \leftarrow 0. 
\]
Hence, the projective dimension of $R/{I_{n,t}}$ equals 2. Since
$\codim(I_{n,t})$ is equal to  the projective dimension of $R/{I_{n,t}}$,
$I_{n,t}$ is Cohen-Macaulay.

\vspace{.2cm}
Let $\Gamma$ be the set of exponents of monomials of $I_{n,t}$. Let 
$\Omega$ be the convex hull of $\Gamma$ in $\mathbb{R}_{\geq 0}^n$. 

\begin{theorem}\label{Thm1}
Let  $(a_1,\dots,a_n) \in \Omega$. Then, for $s \in [n]$ and
$\{i_1,\dots, i_s\} \subseteq [n]$, we have
\[
a_{i_1} + \cdots + a_{i_s} \geq t(s-1).
\]
\end{theorem}

\begin{proof}
Since $\Omega$ is the convex hull of  exponents of monomials of $I_{n,t}$
in $\mathbb{R}_{\geq 0}^n$ and $I_{n,t} = \langle {\bf \tt x}^{e_1},\dots, {\bf \tt x}^{e_n} \rangle$,
there exists, for $i \in [n]$, $c_i \geq 0$ and $d_i \geq 0$ such 
that $c_1+ \cdots + c_n = 1$ and
\[
(a_1,\dots,a_j,\dots , a_n) = c_1e_1 + \cdots + c_ne_n + (d_1,\dots, d_n).
\] 
Thus
\begin{align*}
(a_1,\dots,a_j,\dots , a_n) &\geq c_1e_1 + \cdots + c_ne_n \\
&= (t \sum_{i=2}^n c_i, \dots, t \sum_{\substack{i=1\\i \neq j}}^n c_i, \dots, t \sum_{i=1}^{n-1} c_i).
\end{align*}
In particular, for $j \in [n]$, we have $a_j \geq   t \sum_{\substack{i=1\\i \neq j}}^n c_i$. Therefore, for $\{i_1,\dots, i_s\} \subseteq [n]$, we have
\begin{align*}
a_{i_1}+\cdots+a_{i_s} &\geq t (\sum_{\substack{i=1\\i \neq i_1}}^n c_i+ \cdots +  \sum_{\substack{i=1\\i \neq i_s}}^n c_i)\\
&=  t (-c_{i_1} - \cdots - c_{i_s}+ s \sum_{i=1}^n c_i)\\
&\geq t(s-1).
\end{align*}
\end{proof}

It is well-known, see \cite[Exerice 4.23]{E}, that the integral closure $\overline{I}_{n,t}$ of ${I}_{n,t}$ is the monomial ideal generated by the integral points in $\Omega$.
Thus
\[
\overline{I}_{n,t} = \langle x_1^{a_1}\cdots x_n^{a_n} \, : \, (a_1,\dots,a_n) \in \Omega \cap 
(\mathbb{N} \cup \{0\})^n \rangle.
\]

Let
\begin{align*}
\Delta_{n,t} :=  &\{ (a_1,\dots,a_n) \, : \, a_i\in [t] \mbox{ for all $i$ and } a_1 + \cdots + a_n = t(n-1) \} \\ 
	&\cup \{e_i \, :\, i =1, \dots, n\}.
\end{align*}
The following theorem gives a finite set of generators of $\overline{I}_{n,t}$.

\begin{theorem}\label{main:thm}
Let $\Delta_{n,t}$ and $I_{n,t}$ be as above. Then
\[
\overline{I}_{n,t} = \langle x_1^{a_1}\cdots x_n^{a_n} \, : \, (a_1,\dots,a_n) \in \Delta_{n,t} \rangle.
\]
\end{theorem}

\begin{proof}  
Let $ (a_1,\dots,a_n) \in \Delta_{n,t}$. For $j \in [n]$, let
\[
c_j = 1 - \frac{a_j}{t}. 
\]
It is straightforward to check that $c_j \geq 0$, $c_1+\cdots + c_n = 1$, and
$c_1e_1+\cdots +c_ne_n = (a_1, \dots, a_n)$. Therefore, 
$(a_1,\dots,a_n)$ is a convex combination of $e_1, \dots ,e_n$. 
Thus $(a_1, \dots, a_n) \in \Omega$, and hence 
$x_1^{a_1}\cdots x_n^{a_n} \in \overline{I}_{n,t}$.

To prove the other inclusion,
let $(b_1,\dots, b_n) \in \Omega\cap (\mathbb{N} \cup \{0\})^n $. It is enough to find
$(a_1, \dots, a_n) \in \Delta_{n,t}$ such that
$a_i \leq b_i$ for all $i$. But 
$I_{n,t} \subseteq  \langle x_1^{a_1}\cdots x_n^{a_n} \, : \, (a_1,\dots,a_n) \in \Delta_{n,t} \rangle$,
so we need to prove the above statement only for the case when 
$x_1^{b_1}\cdots x_n^{b_n} \in \overline{I}_{n,t} \setminus I_{n,t}$.
Since $x_1^{b_1}\cdots x_n^{b_n} \notin I_{n,t}$, there exists
$s \geq 2$ and $\{i_1,\dots, i_s\} \subset [t]$ such that
$b_{i_1}, \dots, b_{i_s} < t$. Let $s$ be maximal with this property and
without lose of generality, assume $\{i_1,\dots ,i_s \}=\{1,\dots ,s \}$.
Then $b_1,\dots, b_s < t$ and $b_{s+1}, \dots, b_n \geq t$.
Let $a_{s+1} = ts -(b_1+\cdots+b_s)$. By Theorem \ref{Thm1}, 
$a_{s+1} \leq ts - t(s-1) = t$. Also,
$a_{s+1} > 0$, since $b_{i_1}, \dots, b_{i_s} < t$. Thus 
$a_{s+1} \in [t]$.
It is clear that the $n$-tuple 
$(b_1, \dots, b_s, a_{s+1},t, \dots, t) \in \Delta_{n,t}$ and 
$x_1^{b_1}\cdots x_s^{b_s}x_{s+1}^{a_{s+1}}x_{s+2}^t\cdots x_n^t$ divides
$x_1^{b_1}\cdots x_n^{b_n}$. We conclude that 
$x_1^{b_1}\cdots x_n^{b_n} \in \langle x_1^{a_1}\cdots x_n^{a_n} \, : \, (a_1,\dots,a_n)\in \Delta_{n,t} \rangle.$
\end{proof}

\begin{example}
In Figure \ref{deg4}, the vertices of the triangle correspond to the 
generators of the ideal $I_{3,4}$. Also, all the dots inside and on the boundary
of the triangle correspond to the generators of the integral closure
$\overline{I}_{3,4}$ of $I_{3,4}$.

\vspace{0.1cm}        
\begin{figure}[!htp] 
\centering 
\includegraphics[totalheight=8cm]{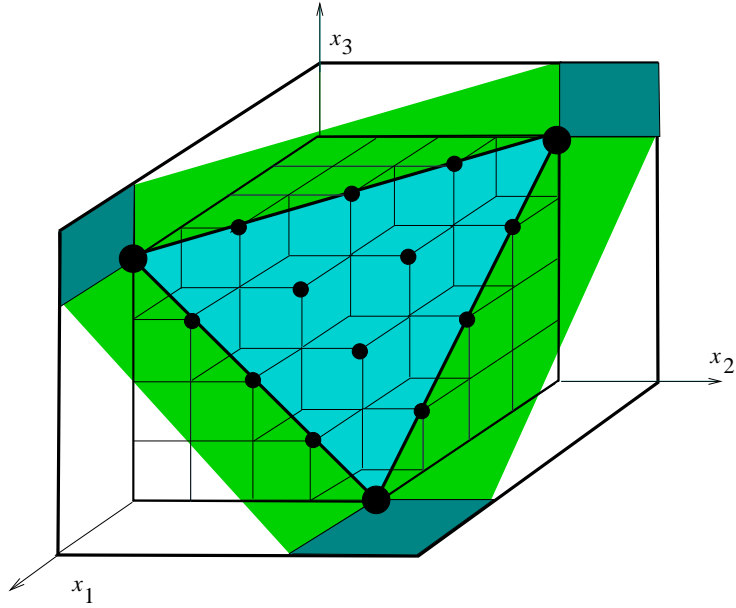} 
\caption{$I_{3,4}$, its convex hull, and the generators of  $\overline{I}_{3,4}$.} 
\label{deg4}
\end{figure}
\end{example}

The following theorem shows that $\overline{I}_{n,t}$ has an embedded associated
prime for $t \geq 2$.

\begin{theorem}\label{int:thm}
For $t \geq 2$, the integral closure $\overline{I}_{n,t}$ has $\langle x_1,x_2,x_3\rangle$
as an embedded associated prime.
\end{theorem}  

\begin{proof}
By Theorem \ref{main:thm}, it is easy to see that
$(\overline{I}_{n,t} \, : \, x_1x_2^{t-1}x_3^{t-1}x_4^t\cdots x_n^t) = \langle x_1,x_2,x_3 \rangle$.
Hence, by \cite[Theorem 4.5]{AM}, 
$\langle x_1,x_2,x_3 \rangle \in \ass(R/\overline{I}_{n,t})$.
Using the same argument, we find that 
$(\overline{I}_{n,t} \> : \> x_2^{t-1}x_3^t\cdots x_n^t) = \langle x_1,x_2 \rangle$ 
and hence 
$\langle x_1,x_2 \rangle \in \ass(R/\overline{I}_{n,t})$. 
Thus $\langle x_1,x_2,x_3\rangle$
is an embedded prime.
\end{proof}

Analogously, $\overline{I}_{n,t}$ has at least $n \choose 3$ embedded
primes, one for each triple of variables. 
Since $\overline{I}_{n,t}$ has embedded primes, by \cite[Theorem 17.3]{Mat},
the following corollary is straightforward.

\begin{corollary}
The ideal $\overline{I}_{n,t}$ is not Cohen-Macaulay for all $t \geq 2$.
\end{corollary}

In the following example, we compute the ideal $\overline{I}_{3,2}$.

\begin{example}
$I_{3,2} = \langle x_1^2x_2^2,x_1^2x_3^2,x_2^2x_3^2 \rangle = 
\langle x_2^2,x_3^2 \rangle \cap  \langle x_1^2,x_3^2 \rangle \cap
 \langle x_1^2,x_2^2 \rangle $. By Theorem \ref{main:thm}, 
it is easy to see that 
\[
\overline{I}_{3,2} = \langle x_1^2x_2^2,x_1^2x_3^2,x_2^2x_3^2, x_1x_2x_3^2,x_1x_2^2x_3,x_1^2x_2 x_3\rangle.
\]
Moreover,
\[
\overline{I}_{3,2}= \langle x_2^2,x_2x_3,x_3^2 \rangle \cap  \langle x_1^2,x_1x_3,x_3^2 \rangle \cap
 \langle x_1^2,x_1x_2,x_2^2 \rangle \cap \langle x_1^2,x_2^2,x_3^2 \rangle.
\]

It is clear that $\langle x_1^2,x_2^2,x_3^2 \rangle$ is an embedded component
of $\overline{I}_{3,2}$.
\end{example}

I have found many other families of examples
but I presented only the simplest one.
In the following example, I  give a Cohen-Macaulay monomial ideal
which is not in the family above and whose integral
closure is not Cohen-Macaulay.

\begin{example}
Let $S = k[x,y,z,w]$ be a polynomial ring over an arbitrary field $k$.
Let $I = \langle x^3yz, x^2w^2,y^2w^3 \rangle =  \langle x^2,y^2 \rangle \cap  \langle x^3,x^2w^2,w^3 \rangle \cap
\langle y,w^2 \rangle \cap  \langle z,w^2 \rangle$. Since 
all the components have codimension 2, 
there are no inclusion relations among them, so by
the uniqueness of minimal primary components, the above 
decomposition of $I$ is irredundant. 
Thus $I$ is a generic (see \cite[Definition 1.1]{MSY}) monomial 
ideal with no embedded primes. By \cite[Theorem 2.5]{MSY}, $I$ is
Cohen-Macaulay. Now by using the computer algebra system 
{\it Normaliz} \cite{BR},
the integral closure $\overline{I}$ of $I$ is
\[
\overline{I} = \langle x^3yz, x^2w^2, y^2w^3, x^2y^2zw, xy^2zw^2, xyw^3\rangle.
\]
It is easy to see that $( \overline{I} \, : \, xyw^2) = \langle x,w \rangle$ and
$( \overline{I} \, : \, xy^2w^2) = \langle x,z,w \rangle$. 
Thus $\overline{I}$ has $\langle x,z,w \rangle$ as an embedded prime and hence 
is not Cohen-Macaulay. 
\end{example}

\section{Acknowledgments}
\noindent
I thank Professor Irena Swanson for suggesting the problem,
and for helpful discussions and conversations. Also, I thank the
referee for pointing out errors in the earlier version of 
this paper.


\end{document}